\documentclass[12pt,a4paper]{article}
\usepackage{amsmath,amssymb,amsfonts}
\def\Bbb{\mathbb}

\title{\bf  Units from square-roots of rational numbers }

\author{Kurt Girstmair}

\date{}

\makeatletter
\let\@@maketitle=\maketitle
\def\maketitle{\def\thispagestyle##1{\relax}\@@maketitle}
\makeatother
\newtheorem{theorem}{Theorem}

\newtheorem{corollary}{Corollary}

\def\BE{\begin{equation}}
\def\EE{\end{equation}}
\def\BD{\begin{displaymath}}
\def\ED{\end{displaymath}}
\def\BA{\begin{array}}
\def\EA{\end{array}}
\def\BEA{\begin{eqnarray*}}
\def\EEA{\end{eqnarray*}}
\def\BI{\bibitem}

\def\N{\Bbb N}
\def\Z{\Bbb Z}

\def\phi{\varphi}
\def\EPS{\varepsilon}

\def\MB{\mbox}
\def\LD{\ldots}

\def\DIV{\,|\,}

\def\MN{\medskip\noindent}
\def\STOP{\hfill$\Box$}

\def\SQDQ{\sqrt{D/Q}}
\def\SQDE{\sqrt{D_1}}
\def\SQD2{\sqrt{D_2}}

\begin{document}

\maketitle

\begin{abstract}

\noindent
Let $D,Q$ be natural numbers, $(D,Q)=1$, such that $D/Q>1$ and $D/Q$ is not a square. Let $q$ be the smallest divisor of $Q$ such that $Q\DIV q^2$. We show that the units $>1$ of the ring
$\Z[\sqrt{Dq^2/Q}]$ are connected with certain convergents of $\SQDQ$. Among these units, the units of $\Z[\sqrt{DQ}]$ play a special role, inasmuch as they correspond
to the convergents of $\SQDQ$ that occur just before the end of each period. We also show that the last-mentioned units allow reading the (periodic) continued fraction expansion of certain
quadratic irrationals from the (finite) continued fraction expansion of certain rational numbers.

\end{abstract}


\section*{1. Introduction and results}

Let $\N=\{1,2,3\LD\}$ denote the set of natural numbers.
Let $D\in\N$, $D$ not a square. In particular, $D>1$. It is well-known how the units in the ring $\Z[\sqrt D]$ can be found, see \cite[p. 93]{Pe}; see also \cite[p. 296]{Ha}, \cite[p. 64]{RoSz}.
Indeed, $\sqrt D$ has a periodic continued fraction expansion
\BE
\label{1.2}
  [b_0,\{b_1,\LD,b_m, 2b_0\}],
\EE
where $b_0,\LD,b_m$ are natural numbers and $\{b_1,\LD,b_m,2b_0\}$ is the period. It is also known that the sequence $b_1,\LD,b_m$ is symmetric, i.e., identical with the sequence $b_m,\LD,b_1$.
Suppose that $m$ has been chosen smallest possible.
Let $r_k/s_k$, $k=0,1,2\LD$ denote the $k$th convergent of $\sqrt D$. As usual, $(r_k,s_k)=1$. Now the units of $\Z[\sqrt D]$ are given by
\BD
\label{1.4}
  r_k+s_k\sqrt D,\enspace k=l(m+1)-1,\enspace l\in\N.
\ED
In other words, the convergents occurring just before the end of each period supply these units.

In this note we study numbers $D/Q$, where $D, Q\in\N$ are such that $(D,Q)=1$, $D/Q$ is not a square, and $D/Q>1$, i.e., $Q<D$. The continued fraction expansion of $\sqrt{D/Q}$ also has the form
(\ref{1.2}), and $b_1,\LD,b_m$ have the same symmetry property. The following theorem describes the units obtained from convergents of the continued fraction expansion (\ref{1.2}) of $\SQDQ$.
To this end let $q$ be the smallest divisor of $Q$ such that $Q\DIV q^2$. Then $D_1=Dq^2/Q$ is a natural number $\ge D/Q$ and not a square. For the time being, we do not assume that $m$ is smallest possible.

\begin{theorem} 
\label{t1}

In the above setting, let $r+s\SQDE$ be a unit in $\Z[\SQDE]$. Let $t=(r,q)$. Then
$(r/t,sq/t)$ is a convergent of $\SQDQ$.
Moreover, if $(r/t,sq/d)=(r_k,s_k)$, then $(r_{k+m+1},s_{k+m+1})$ also has the form $(r'/t,s'q/t)$ for another unit $r'+s'\SQDE$ of $\Z[\SQDE]$.
Here $(r',q)=t$ for the same number $t$.

\end{theorem} 

\MN
So the units of $\Z[\SQDE]$ are closely connected with certain convergents of $\SQDQ$.
Let $D_2=DQ$. Since $D_2=D_1(Q/q)^2$, the ring $\Z[\SQD2]$ is contained in $\Z[\SQDE]$. Its units play a special role in the continued fraction expansion of $\SQDQ$.

\begin{theorem} 
\label{t2}

In the above setting, let $m$ be smallest possible.
Then the units of $\Z[\SQD2]$ have the form
$r+s\SQDE$,
where $r=r_k, sq=s_k$, $k=l(m+1)-1, l\in\N$; i.e., the numbers $r_k/s_k$ are the convergents of $\SQDQ$
occurring just before the end of each period of this quadratic irrational.
\end{theorem} 

\MN
The above theorems suggest the following notation. A unit of $\Z[\SQDE]$ that lies in $\Z[\SQD2]$ is called a {\em regular} unit, whereas a unit of $\Z[\SQDE]$ that is not in $\Z[\SQD2]$
is called an {\em irregular unit}.

Irregular units are only possible if $q^2\ne Q$. Moreover, they can only occur if the unit group of $\Z[\SQDE]$ is larger than that of $\Z[\SQD2]$.

\MN
{\em Example.} Let $D=157$, a prime, and $Q=45$. Then $q=15$. Here $D_1=5\cdot 157=785$ and $D_2=9\cdot 785$. So we have $\Z[\SQDE]=\Z[\sqrt{785}]\neq\Z[\SQD2]=\Z[3\sqrt{785}]$.
The fundamental unit of $\Z[\SQDE]$ is $\EPS=28+\sqrt{785}$. Since $\EPS\not\in\Z[\SQD2]$, $\EPS$ is an irregular unit.

Now $\SQDQ=\sqrt{157/45}$ has the smallest possible period length $16$ and $r_{15}=4923521$, $s_{15}=2635920$. Here $\EPS^4=4923521 + 175728\sqrt{785}=r+s\SQDE$, and $r_{15}=r$, $s_{15}=sq=175728\cdot 15$.
Theorem \ref{t2} says that $r+s\SQDE$ is a regular unit. Indeed, we have
$r+s\SQDE=4923521+58576\SQD2$. Of course, $r+s\SQDE$ is the fundamental unit of $\Z[\SQD2]$.

We have $r_4=28$ and $s_4=15$, whose connection with $\EPS$ is obvious. We obtain $\EPS^2=1569+56\SQDE$. Here $t=(1569, q)=3$. This gives $r_7=1569/3=523$ and $s_7=56\cdot 15/3=280$. Finally,
$\EPS^3=87892 + 3137\SQDE$ and $r_{10}=87892, s_{10}=47055=3137\cdot 15$. Hence we have irregular units connected with the convergents of indices $2,7,10$ and $2+l\cdot 16$, $7+l\cdot16$, $10+l\cdot 16$, $l\in\N$.
The regular units are connected with the convergents of indices $15$ and $15+l\cdot 16$, $l\in\N$.

\MN
Of course, the regular units can also be found by the continued fraction expansion of $\SQD2$. Observing the monotonous growth of the numerators and denominators of convergents, we obtain
the following corollary.

\begin{corollary}
\label{c1}
In the above setting, let $m$ be smallest possible such that $\SQDQ=[b_0,\{b_1,$ $\LD, b_m,2b_0\}]$. Let $n$ be smallest possible such that  $\SQD2=[b'_0,\{b'_1,\LD,b'_n,2b'_0\}]$.
Let $r_k/s_k$ be the $k$th convergent of $\SQDQ$ and $r'_k/s'_k$ the $k$th convergent of $\SQD2$. Then
\BD
 r_{l(m+1)-1}=r'_{l(n+1)-1},\enspace s_{l(m+1)-1}=s'_{l(n+1)-1}Q
\ED
for all $l\in\N$.

\end{corollary}

\MN
In the above example, $m=16$ and $\SQD2$ has the smallest possible period length $n=8$. So we have $r'_7=r_{15}=4923521$ and $s'_7=58576$, $s'_7Q=s_{15}=2635920$.

Theorems \ref{t1} and \ref{t2} exhibit the structure of the units $>1$ in $\Z[\SQDE]$. To this end let $m$ be smallest possible.
If $k=m$, then the convergent $r_k/s_k$ of $\SQDQ$ defines the fundamental unit $\eta$ of $\Z[\SQD2]$, by Theorem \ref{t2}. Suppose that the fundamental unit $\EPS$ of $\SQDE$ is different from $\eta$.
Then $\eta=\EPS^j$ for some $j>1$. By Theorem \ref{t1}, the units $\EPS, \EPS^2,\LD,\EPS^{j-1}$ are connected with convergents $r_{k_1}/s_{k_1}, \LD, r_{k_{j-1}}/s_{k_{j-1}}$ of $\SQDQ$,
and, due to the strictly monotonic growth of $r_k+s_k\SQDQ$ for $k\to\infty$, we have $k_1<\LD<k_{j-1}<m$. By Theorem \ref{t1}, the convergents $r_{k_{1+m+1}}/s_{k_{1+m+1}}$, \LD, $r_{k_{j-1+m+1}}/s_{k_{j-1+m+1}}$
define units in $\Z[\SQDE]$ larger than $\eta$ but smaller than $\eta^2$, the latter being defined by $r_{2m+1}/s_{2m+1}$. Hence there is only one possibility, namely, the aforesaid convergents define the
units $\EPS^{j+1},\LD,\EPS^{2j-1}$, which are the only units of $\Z[\SQDE]$ that are larger than $\eta$ but smaller than $\eta^2$. In the same way one can proceed for higher powers of $\EPS$.

In what follows we also consider the continued fraction expansion of a {\em rational} number $x$,
so
\BE
\label{1.5}
 x=[c_0,\LD,c_n].
\EE
Here $c_n$ need not be $\ge 2$. Depending on the situation, we require, instead, that $n$ is even or odd. This can be obtained on replacing $c_n$ by the sequence $c_n-1,1$ at the end of (\ref{1.5}).

\begin{theorem} 
\label{t3}
Let $r_k+s_k\SQDQ$ be a regular unit, i.e., $k=l(m+1)-1$ if $m$ is smallest possible.  Let $t$ be a divisor of $s_k$.
Let
\BD
\label{1.8}
  r_k/t=[c_0,\LD,c_n],
\ED
where $n\equiv k\mod 2$.
Then
\BD
\label{1.10}
  (s_k/t)\SQDQ=[c_0,\{c_1,\LD,c_n, 2c_0\}].
\ED
\end{theorem} 

\MN
In other words, we may read the continued fraction expansion of the quadratic irrational $(s_k/t)\SQDQ$ from the continued fraction expansion of the rational number $r_k/t$.

\MN
{\em Example.} In the above example, we have $r_{15}=4923521$, $s_{15}=2635920=16\cdot 9\cdot 5\cdot 7\cdot 523$. We choose $t=16\cdot 9\cdot 7$ and obtain
$r_{15}/t=[4884,2,4,12,4,2]$. This yields
\BD
  (s_{15}/t)\sqrt{D/Q}=523\sqrt{785}/3=[4884,\{2,4,12,4,2,9768\}].
\ED

\section*{2. Proofs} 

{\em Proof of Theorem 1.} The number $r+s\SQDE$ is a unit, if, and only if,
\BE
\label{2.2}
 r^2-s^2D_1=\pm 1.
\EE
Since $D_1=Dq^2/Q$ and $t=(r,q)$, this is the same as saying
\BE
\label{2.4}
  (r/t)^2-(sq/t)^2D/Q=\pm(1/t)^2.
\EE
Since $1/t\le 1$ and $D/Q>1$, we see that $r/t=r_k$, $sq/t=s_k$ defines a convergent of $\SQDQ$, see \cite[p. 39]{Pe}.

For the next step we need the $k$th complete quotient $x_k$ of $x_0=\SQDQ$. By \cite[p. 80]{Pe},
\BE
\label{2.5}
  x_k=(\SQD2+P_k)/Q_k
\EE
with $P_k\in\Z$ and $Q_k\in \N$. In particular, $x_0=\SQDQ=(\SQD2+P_0)/Q_0$, with $P_0=0$ and $Q_0=Q$.
From \cite[p. 71]{Pe} we obtain
\BE
\label{2.6}
  r_k^2-s_k^2D/Q=\pm Q_{k+1}/Q.
\EE
Suppose now that $r_k=r/t$, $s_k=sq/t$. Comparing (\ref{2.4}) with (\ref{2.6}), we see that
\BD
  Q_{k+1}/Q=1/t^2, \MB{ i.e. }, Q_{k+1}=Q/t^2.
\ED
In particular, $t^2\DIV Q$. Now the complete quotient $x_{k+1+m+1}$ is the same as $x_{k+1}$, whence we obtain
\BD
  r_{k+m+1}^2-s_{k+m+1}^2D/Q=\pm Q_{k+1}/Q=\pm(1/t)^2.
\ED
This gives
\BD
  (r_{k+m+1}t)^2-(s_{k+m+1}t)^2D/Q=\pm 1.
\ED
Therefore, however, $Q$ divides $(s_{k+m+1}t)^2$, which, by the minimal property of $q$, implies $q\DIV s_{k+m+1}t$. Accordingly, we may define $r',s'\in\Z$
by $r'=r_{k+m+1}t$, $s'=s_{k+m+1}t/q$ and obtain
\BD
   r'^2-s'^2Dq^2/Q=r'^2-s'^2D_1=\pm 1.
\ED
This is the desired result.
\STOP

\MN
{\em Proof of Theorem \ref{t2}.}

Let $r+s\SQDE$ be  a unit in $\Z[\SQDE]$. It is easy to see that this unit lies in $\Z[\SQD2]$ if,
and only if, $Q\DIV sq$.

Suppose that this holds. Then the identity (\ref{2.2})
implies $(r,Q)=1$. By the argument used in the proof of Theorem \ref{t1}, $r_k=r, s_k=sq$ is a convergent of $\SQDQ$.
From (\ref{2.6}) we obtain
\BD
  r_k^2-s_k^2D/Q=\pm Q_{k+1}/Q=\pm 1.
\ED
Therefore, the complete quotient $x_{k+1}$ has the form
\BE
\label{2.8}
  x_{k+1}=(\sqrt{D_2}+P)/Q,
\EE
for some integer $P$ (see (\ref{2.5})).

We show that $Q\DIV P$.
Since $x_0=\SQDQ=[b_0,\LD,b_k,x_{k+1}]$, we obtain
\BD
   x_0=\frac{r_k x_{k+1}+r_{k-1}}{s_k x_{k+1}+s_{k-1}}.
\ED
Let
\BD
\left(
      \begin{array}{cc}
        a & b \\
        c & d \\
      \end{array}
    \right)
                        =\left(
                               \begin{array}{cc}
                                 r_k & r_{k-1} \\
                                 s_k & s_{k-1} \\
                               \end{array}
                             \right),
\ED
a matrix whose determinant $ad-bc=\pm 1$. In particular, $c= s_k$. We use $x_0=\SQD2/Q$ and (\ref{2.8}) in order to compare the coefficients in the relation
\BD
   x_0(cx_{k+1}+d)=ax_{k+1}+b.
\ED
This yields
\BE
\label{2.10}
  b=-aP/Q+cD_2/Q^2=-aP/Q+cD/Q.
\EE
Now we know $Q\DIV c$. Then $cD/Q\in\Z$ and, by (\ref{2.10}), $-aP/Q\in\Z$. But $(a,Q)=1$, since, otherwise, we have a contradiction to $ad-bc=\pm 1$.
Therefore, $Q\DIV P$ and $x_{k+1}=x_0+C$, $C\in\Z$. But then $x_{k+2}=x_1$, which means that
$k+1$ marks the end of a period. Thus, $k=l(m+1)-1$, $l\in\N$.

Conversely, let $k=l(m+1)-1$, $l\in\N$. Since $x_{k+2}=x_1$, we have $x_{k+1}=x_0+C$, $C\in\Z$. On the other hand, $x_{k+1}$ has the form
of (\ref{2.8}), by \cite[p. 81]{Pe}. This yields $Q\DIV P$.
Then $-aP/Q\in\Z$ and, by (\ref{2.10}), $cD/Q\in\Z$.  Since $(D,Q)=1$,
$Q\DIV c$, and $c=s_k$. Again, we have (\ref{2.6}) with $Q_{k+1}=Q$.
Accordingly, $r_k/s_k$ belongs to a unit in $\Z[\SQD2]$.
\STOP

\MN
{\em Proof of Theorem \ref{t3}.}
The regular unit $r_k+s_k\SQDQ$ can be written
\BD
  r_k+s_k\SQDQ=r_k+t\sqrt{(s_k/t)^2D/Q}.
\ED
If we write $(s_k/t)^2D/Q$ in the form $D'/Q'$, $(D',Q')=1$, then
\BE
\label{2.12}
   t\sqrt{D'/Q'}=a\sqrt{D'Q'}
\EE
for some natural number $a$. This identity will be proved below. It says that $r_k+t\sqrt{D'/Q'}$ is a regular unit (with respect to $D'$ and $Q'$).
Accordingly, Theorem \ref{t2} yields
\BD
\label{2.14}
  r_k/t=[c_0,\LD,c_n],
\ED
where
\BD
\label{2.16}
  (s_k/t)\SQDQ=\sqrt{D'/Q'}=[c_0,\{c_1,\LD,c_n, 2c_0\}]
\ED
($n$ not necessarily smallest possible). Here $n\equiv k\mod 2$ since the sign of $r_k^2-t^2D'/Q'=r_k^2-s_k^2D/Q$ ($=\pm 1$) defines the parity of $k$ and $n$.

In order to prove (\ref{2.12}), we consider the $p$-exponent $v_p(a)$ of a rational number $a$. Because $Q'\DIV Q$, only prime divisors $p$ of $Q$ are relevant. Let $v_p(Q)=e_p$, $v_p(t)=f_p$, and $v_p(s_k/t)=g_p$.
Since $Q\DIV s_k$, we have $f_p+g_p\ge e_p$.
Now $v_p((s_k/t^2)D/Q)=2g_p-e_p=v_p(D'/Q')$. Hence
\BD
  v_p(D')=\left\{
            \begin{array}{ll}
              2g_p-e_p, & \MB{ if } e_p\le 2g_p; \\
              0, & \MB{ if } e_p>2g_p.
            \end{array}
          \right.
\ED
Further,
\BD
  v_p(Q')=\left\{
            \begin{array}{ll}
              0, & \MB{ if } e_p\le 2g_p; \\
              e_p-2g_p, & \MB{ if } e_p>2g_p.
            \end{array}
          \right.
\ED
Accordingly,
\BE
\label{2.20}
  v_P(D'Q'))=\left\{
            \begin{array}{ll}
              2g_p-e_p, & \MB{ if } e_p\le 2g_p; \\
              e_p-2g_p, & \MB{ if } e_p>2g_p.
            \end{array}
          \right.
\EE
Now
\BE
\label{2.22}
  v_p(t^2D'/Q')=v_p(t^2(s_k/t)^2D/Q)=2f_p+2g_p-e_p.
\EE
We obtain, from (\ref{2.20}) and (\ref{2.22}),
\BD
 v_p\left(\frac{t^2D'/Q'}{D'Q'}\right)=\left\{
            \begin{array}{ll}
              2f_p, & \MB{ if } e_p\le 2g_p; \\
              2f_p+4g_p-2e_p, & \MB{ if } e_p>2g_p.
            \end{array}
          \right.
\ED
Because $f_p+g_p\ge e_p$, this exponent is always nonnegative and even. This implies (\ref{2.12}).
\STOP


\vspace{1cm}
\noindent
Institut f\"ur Mathematik \\
Universit\"at Innsbruck   \\
Technikerstr. 13/7        \\
A-6020 Innsbruck\\
Kurt.Girstmair@uibk.ac.at

\end{document}